\documentclass[11pt,letterpaper]{article}
\usepackage[T1]{fontenc}
\usepackage[utf8]{inputenc}
\usepackage{mathpazo}
\usepackage[margin=1in]{geometry}
\usepackage{microtype}
\usepackage{parskip}
\usepackage[hidelinks]{hyperref}
\usepackage{bookmark}
\ifdefined\pdfgentounicode
  \input{glyphtounicode}
\fi
\hypersetup{
  pdftitle={University Research and Leadership in the Age of AI},
  pdfsubject={Mobilizing Intellectual Capacity and Building the Infrastructure for Discovery},
  pdfauthor={Bacim Alali},
  bookmarksnumbered=true
}
\title{\Large\bfseries University Research and Leadership\\
in the Age of AI\\[0.5em]
\large\normalfont Mobilizing Intellectual Capacity and Building\\
the Infrastructure for Discovery}
\author{Bacim Alali\\[0.3em]
{\small\href{mailto:bacimalali@math.ksu.edu}{bacimalali@math.ksu.edu}}}
\date{}

\begin{document}
\maketitle

\begin{abstract}
\noindent
Artificial intelligence (AI) is expected to disrupt research across academia, changing how knowledge is produced and how research is organized. Recent developments in mathematics provide an early example of strong AI research capabilities and the professional disruption accompanying them. Universities should respond by building the scientific, computational, and organizational capacity to develop AI and shape the research it enables. The relationship is reciprocal: AI can advance disciplinary research, while scientific questions, experimental evidence, and disciplinary expertise can guide AI development. This article proposes connecting shared computational resources with sustained technical teams, interdisciplinary research, experimental facilities, and knowledge infrastructure that makes findings traceable and reusable. It distinguishes three commitments with different resource requirements: operating existing models, developing AI for disciplinary research, and investigating the foundations of AI. Commercial services can support this agenda alongside institutional capabilities for independent investigation, evaluation, and preservation of the research record. Consortia and national facilities can broaden participation and support capabilities beyond individual campuses. Initial programs should assess scientific value, total costs, and contributions to AI development, with expansion guided by evidence. These changes require sustained funding, professional support, and recognition of the work that makes research reliable and cumulative. Building this capacity will enable universities to pursue new research directions, train future researchers, and exercise scientific leadership as AI transforms academic inquiry.
\end{abstract}

\section{Introduction}

Artificial intelligence is expected to disrupt research across academia. Its growing capabilities are likely to change the questions researchers pursue, the ways they conduct and organize their work, and the computational resources and expertise they need. Universities therefore face decisions about how to support their researchers, develop new scientific capabilities, and participate in shaping the technology on which their work increasingly depends. The timing and form of this disruption will differ across disciplines, but developments in fields already experiencing it can help clarify what universities should prepare for.

I approach these questions as a mathematician. Mathematics provides a current example of both the strength of AI's research capabilities and the disruption accompanying their development. Recent work documents AI systems producing proofs and counterexamples and making substantial progress on long-standing open problems across several mathematical areas~\cite{aletheia,openai-ten,anthropic-zeta,openai-math-release}. These developments are already prompting mathematicians to reconsider research directions, standards of contribution, and the future of their profession. The Advisory Group on Mathematics and Artificial Intelligence has called for release practices that preserve human understanding, verification, and responsibility for mathematical arguments~\cite{agmai-release}. Students and early-career researchers face pressing questions about their careers, while researchers using AI are finding opportunities to pursue a wider range of mathematical investigations.

The experience of mathematics makes a broader institutional challenge concrete. Taking advantage of powerful AI requires dependable access to computation and suitable models, technical support, and resources for examining, preserving, and building on research results. Universities should help develop this infrastructure and the AI methods it supports. Their disciplinary expertise and experimental facilities also provide opportunities to shape AI development: scientific questions can motivate new methods, and evidence from research can reveal how those methods should improve. The relationship is reciprocal---AI can advance disciplinary research, and disciplinary research can advance AI. This paper uses developments in mathematics as a starting point for arguing that universities should build the capacity to lead this work across their research communities.

The same institutional questions take different forms across disciplines. Mathematical arguments can be assessed through proof, while claims about physical or biological systems require observational or experimental evidence. Research in the humanities and social sciences brings further requirements concerning sources, interpretation, and context. These differences should shape the tools and standards universities develop. The mathematical example motivates preparation; the value of a particular AI system in another field must be established through that field's research practices.

Building this capacity takes years of sustained scientific and technical work. Universities can undertake it through consortia, national facilities, and partnerships, with commitments suited to their research strengths. The focus here is universities as research institutions. Training future researchers is part of that capacity; a broader redesign of curricula lies outside this paper's scope. The aim is to enable universities to initiate research, develop its instruments, and help determine its direction as AI capabilities expand.

\section{Mobilizing intellectual capacity}

The distinction between using AI and developing AI is central to this argument. An \emph{AI model} learns patterns from data or examples to produce outputs such as predictions, text, or proposed designs. It can also serve as a scientific model, for example by predicting a material's properties. An \emph{AI research system} combines such models with the software, information sources, scientific tools, and instruments through which researchers use them. Universities should contribute to the development of both models and research systems.

The case for that contribution begins with the breadth of university research. Experimental scientists design tests that distinguish competing explanations. Engineers determine whether a design meets physical, manufacturing, and safety constraints. Biomedical researchers investigate why an apparent relationship holds in one setting and fails in another. Scholars in the humanities and social sciences examine interpretation, institutions, incentives, and the purposes to which technologies are directed. These researchers should help define the problems AI developers work on and the standards by which progress is judged, rather than be consulted after those decisions have been made.

The relationship works in both directions. Consider a materials project in which AI helps researchers select compounds to test. If its predictions repeatedly fail for a particular class of materials, experiments may reveal a limitation in the information the model uses or the way it learns physical relationships. Materials researchers and AI researchers can then build a better model and test whether it predicts previously unexamined materials more accurately. In engineering, an AI-generated design that keeps violating manufacturing constraints motivates methods that incorporate those constraints during generation rather than catching every violation afterward. These are research programs in AI development, not additional uses of an unchanged product.

A historian working with a multilingual archive might help develop AI methods that retain source attribution and distinguish documentary evidence from interpretation. Difficulties with incomplete records or changing terminology would then become research questions for both historians and AI researchers.

These collaborations require shared questions, sustained resources, and responsibility for the resulting work. A \emph{benchmark}, a standardized set of tasks used to compare AI performance, can inform evaluation, but scientific value also depends on what the work establishes and enables. A useful new method, an informative negative result, or a well-supported explanation may matter more than a gain on a general ranking. An AI system that makes a previously impractical environmental study possible should be assessed through the quality and value of that study.

No small group of organizations can represent the full range of worthwhile questions. AI companies produce fundamental discoveries and substantial public benefit, but their priorities cannot substitute for those of the wider research community. Diseases that receive little commercial investment, local water and environmental problems, and engineering needs in poorly resourced settings all deserve attention. University participation should broaden the people and purposes shaping AI, not replace corporate concentration with concentration in a few wealthy campuses.

\section{AI as research infrastructure}

Computational AI resources should become a dependable part of the research environment, alongside libraries, scientific instruments, and research computing. Researchers should have a baseline of suitable AI models and technical support, free at the point of use, rather than depend on personal subscriptions. That baseline should include interactive assistance, access from research software, and environments in which AI can carry out longer sequences of work using approved tools. A researcher might connect an AI assistant to a simulation and ask it to compare proposed designs within a defined computing allowance. Staff must be available to configure such work, diagnose failures, and protect research data; an account alone does not provide that capacity.

For many institutions, a consortium or state university system is the practical unit of organization. Shared facilities can pool demand and sustain engineering teams beyond the means of individual departments. Empire AI brings together New York universities and research institutions to share advanced computing for AI research in the public interest~\cite{empireai}. Its development illustrates how a consortium can provide capabilities beyond those readily available to individual campuses. The scientific value of such investment depends on the research teams, technical support, and collaborative programs that use it.

The U.S. National Artificial Intelligence Research Resource (NAIRR) connects research communities to computing, AI models, data, software, and expertise~\cite{nairr}. The Department of Energy's Genesis Mission aims to connect advanced computing, AI models, scientific data, and experimental facilities in a coordinated environment for discovery; its consortium includes academic institutions and identifies AI model development and validation, data standards, computing infrastructure, and laboratory automation as areas for joint work~\cite{genesis,genesisconsortium}. Both are institutional commitments to integration; neither is yet evidence that their scientific ambitions have been met. Participation should complement sustained academic capacity: universities need their own expertise and scientific direction to shape joint research, pursue investigator-led questions, and evaluate results independently.

Academic control over part of this infrastructure matters. During training, an AI model's numerical weights are adjusted using data or examples. When these weights are released, researchers can run the model themselves, preserve a particular version, examine its internal behavior, and modify it through further training, subject to its license. This is commonly called an \emph{open-weight} model. Such access supports substantial research, including experiments that a commercial interface accepting questions and returning answers may not permit.

The release of weights alone does not disclose how a model was developed. Studying the influence of training data or comparing training procedures can require access to datasets, training code, documentation, and saved model versions from different stages of development. The Model Openness Framework~\cite{mof} distinguishes levels of openness according to the materials released and the permissions for their use. Universities should assess what a release makes available against the research questions they intend to investigate.

The appropriate arrangement will usually combine university-operated services, rented computation, national facilities, and commercial AI. Access arrangements must also support the long-term preservation of research records, a responsibility considered in Section~4.

\section{The case for buying access, and its limits}

A serious objection to academic AI infrastructure concerns opportunity cost. Commercial services can provide capable AI without a university bearing the full cost of developing and operating the underlying models. Shared purchasing can broaden access while leaving resources available for instruments, research staff, and graduate support. The concentration of computation in industry strengthens this economic objection: academic teams face constraints in compute-intensive AI research, while commercial providers can spread development and operational costs across many users~\cite{computedivide,industryai,naturecompute}. A proposal for academic operation should therefore compare purchased services, shared facilities, and local systems against the intended workloads, including technical support, security, continuity, integration, and maintenance.

For many routine activities, purchased access will be a practical choice. The appropriate arrangement depends on the institution's research requirements and existing resources. Service procurement and research investment nevertheless answer different questions. Buying useful tools supports current work; developing AI creates capabilities, knowledge, and methods that current products may not provide. Universities have reasons to pursue that research even when commercial services serve their everyday needs well. Three continuing responsibilities make the limits of access alone particularly clear.

First, universities should establish the research access they need before committing to a service or partnership. Agreements should specify which materials researchers can examine and which experiments or modifications they can undertake.

Second, universities have a continuing responsibility for the research record. A published mathematical proof can remain independently checkable after the model that discovered it has been retired. Reproducing a claim about that model's discovery process may additionally require retained instructions, retrieved information, software, model access, and execution conditions. Similarly, an experimentally established property of an AI-proposed material can remain valid when its generating system changes. Preservation arrangements should specify which records and checks remain available, who maintains them, and what can still be examined when a service changes. These obligations require durable institutional commitments, including appropriate agreements with commercial providers.

Third, independent evaluation requires appropriate access and authority to publish findings, including unfavorable ones. Purchased services, open research artifacts, and negotiated access can each support such work under suitable conditions. Section~9 examines the institutional arrangements needed to sustain that independence.

These responsibilities help define the research access, technical expertise, and preservation capacity that universities need alongside ordinary services. The positive research agenda is broader: disciplinary and foundational AI research can create new capabilities and open new scientific questions. Its scale should follow the questions being investigated and the resources required to answer them. Consortia, national facilities, and industry partnerships offer different ways to support that work.

\section{Developing AI, not only deploying it}

Three commitments with very different costs and staffing requirements give universities a structure for deciding what to attempt: operating existing AI models, developing AI for particular research needs, and investigating the foundations of AI. They need not be undertaken at the same scale or by the same institutions.

The first is \emph{shared operation of existing AI models}. Running suitably licensed models on academically controlled infrastructure can support controlled experiments, continued access to selected versions, and appropriate handling of sensitive information. Such services require maintained software and execution environments as well as model weights. Their scope should follow demonstrated research needs and comparisons with purchased services. A shared service should evaluate models on relevant tasks and revise its offerings as capabilities, costs, and requirements change.

The second is \emph{disciplinary AI development}. Teams can adapt an existing model by further training, build a specialized model, or connect AI to scientific references, simulations, and instruments. These are different activities: giving a model access to a reliable database may improve its answers without changing its weights. Scientific AI need not consist only of language models. Other AI models predict physical properties, analyze combinations of images and measurements, or propose a material or device with desired characteristics. Universities should establish which approaches work for particular scientific questions and which limitations demand new methods.

The third is \emph{foundational AI research}: investigating how AI models are organized, how they represent information, how they learn, why they fail, and how their behavior can be understood and directed. Some questions can be studied with modest computation; others require substantial shared facilities. Training a model from scratch means developing its learned settings without starting from an existing model's weights. University teams can undertake such work at scales suited to their questions. Repeatedly training very large general-purpose models requires a much larger commitment of people, equipment, and experimental capacity. National centers and major collaborations can support that ambition where the scientific case warrants it.

Existing projects illustrate different parts of this agenda. Apertus~\cite{apertus} is a language model with openly released weights, source code, and documentation of training methods and datasets. It demonstrates that universities working with a national computing center can develop the underlying technology. Olmo~3, from the nonprofit Allen Institute for AI, releases model weights together with training data, training code, and saved model versions from different stages of development~\cite{olmo3}. These materials enable researchers to investigate how capabilities emerge during training and to test changes to the development process.

Marin began at Stanford University and is now primarily developed by Open Athena with community contributors. It shares model weights, code, training-data sources, and records of successful and unsuccessful experiments, making its development process available for others to examine and extend~\cite{marin}.

These commitments are not a ladder every institution must climb. A university team developing a new learning method may use national computing facilities without operating a local AI service. A consortium may provide broad access while supporting selected disciplinary projects. Across these arrangements, long-term plans should account for changes in what external providers release and sustain the expertise needed to maintain and develop the models on which research depends.

\section{Integrating AI across the research process}

Investment should connect the stages of work within a research project. Researchers can investigate how AI contributes at each stage of inquiry: interpreting prior work, formulating questions, developing explanations, designing experiments, performing calculations, checking results, and communicating findings. Connecting these activities requires keeping track of the evidence and assumptions used at each stage. If an AI-assisted analysis relies on a particular dataset or assumption, that information must remain available when the analysis is used to plan an experiment or support a conclusion.

An AI system proposing an experiment should report the uncertainty behind its prediction. A simulated result should remain distinguishable from an actual measurement. An engineering recommendation should identify the conditions under which it is expected to work. If an instrument is later found to be miscalibrated, researchers should be able to trace which analyses and conclusions used its measurements. These requirements connect the design of research software to the evidence on which scientific judgments depend.

The combined research environment should itself be tested. Does AI assistance help choose more informative experiments? Does a model use a simulator correctly? Does automated checking find important errors without generating so many false alarms that it costs more time than it saves? Evaluation must include the work needed to check and repair AI output. A larger volume of proposed explanations, untested designs, or preliminary analyses is not a larger body of reliable knowledge.

Connecting AI to research also creates opportunities to improve it. Failed predictions or designs can reveal weaknesses in a model, gaps in the information available to it, or problems in its use of scientific software. Investigating these failures can produce better research tools and new AI methods.

\section{Laboratories as generators of new evidence}

University laboratories allow researchers to ask questions that existing data do not answer and generate evidence to investigate them. In biotechnology, nanotechnology, physics, materials science, and engineering, this capacity includes instruments and physical space, and also technical staff, reliable procedures, and expertise in interpreting measurements. These resources connect AI's suggestions about the world to tests in the world.

CRESt --- Copilot for Real-world Experimental Scientists --- provides a published example of this integration~\cite{crest}. It combines models that process multiple kinds of scientific information, methods for selecting promising experiments, and robotic testing for catalyst discovery. Lila Sciences~\cite{lila,mitlila} illustrates commercial investment in AI and automated laboratories. Extending such approaches to other laboratories requires adaptation to their instruments, workflows, and scientific questions, with evidence of the resulting benefits.

Universities should choose the scientific questions, determine which experiments provide informative evidence, and develop AI methods that use that evidence. Such work can produce new findings, better scientific explanations, and improved methods for selecting experiments. Partnerships can contribute resources and expertise, but agreements must preserve academic researchers' ability to investigate how the resulting AI performs and to publish what they find. Supplying experiments while surrendering influence over their research use is not leadership.

A university program need not begin with an autonomous laboratory. It can connect selected instruments, improve experimental records, support AI-assisted design, and automate well-defined procedures under oversight. Extending such work requires records that different laboratories can exchange, calibrated equipment, maintenance, safety controls, and staff able to make the connections reliable. Existing facilities are an advantage, not a cost-free robotic research network. Permission for AI to propose an experiment must remain distinct from permission to operate equipment or initiate a procedure.

Faculty, research scientists, postdoctoral researchers, graduate and undergraduate researchers, and technical specialists can all contribute. Designing informative experiments and producing carefully documented datasets for AI development are scientific contributions deserving recognition, not a supply of inexpensive data-generation labor. Stable technical staff are essential: the ability to maintain instruments and interpret their records must survive the end of an individual project.

The value of experimental data depends on the questions they help answer. Measurements can be biased, redundant, or poorly documented. A small set of well-chosen experiments can be more informative than a large volume of measurements that do not distinguish competing explanations.

\section{Knowledge infrastructure for cumulative discovery}

The records produced within individual projects should become resources that other researchers can examine and build on. This requires shared knowledge infrastructure that connects findings to their assumptions, evidence, methods, and subsequent corrections, while preserving attribution. Papers remain essential sources and explanations; explicit relationships can make their contributions easier to use in new investigations.

Mathematics provides a concrete example. Researchers often need to determine whether an estimate applies under their hypotheses, whether two formulations express the same result, or whether an alternative proof exposes a transferable method. Recovering these relationships repeatedly from separate texts consumes effort and can conceal useful connections. A companion proposal describes an infrastructure for linking mathematical statements, proof routes, and methods to their original sources~\cite{alaliinfrastructure}. Its value should be evaluated against strong literature-search and retrieval alternatives. The institutional aim is to make research cumulative and reusable, with the technical design revised in light of what works.

In an experimental field, a usable research record connects a conclusion to its protocol, measurement conditions, uncertainty, data processing, and unsuccessful as well as successful tests. A model asked to recommend an experiment needs to distinguish a result obtained under narrowly controlled conditions from one supported across many settings. Researchers evaluating an AI-generated analysis need access to the evidence it used and the transformations applied to it. Such records can also supply well-documented material for developing and evaluating AI, subject to appropriate rights and independent assessment.

Formal proof assistants illustrate the lasting value of shared research tools. Lean originated at Microsoft Research, and its development and mathematical libraries have drawn contributions from academic, industrial, and nonprofit researchers~\cite{lean-history}. A formal checker can establish that a particular statement follows from declared definitions and assumptions; researchers must also assess whether the statement faithfully captures the intended mathematics. In empirical research, formal checks can validate aspects of an analysis, while its observational or experimental basis requires separate assessment. Across these settings, the common requirement is an inspectable connection between a claim and its supporting evidence.

AI can help build these records by extracting information and proposing links between findings. Each relationship should retain its sources and indicate what has been checked. Software may suggest that two studies used comparable measurements when a change in procedure makes their results incompatible. Corrections should identify the analyses or proof routes that require review, while preserving conclusions supported by independent evidence. These functions remain useful as AI improves at verification and explanation: their purpose is to make the products of research available across projects and changing tools.

Making research reusable also requires clear conditions of use. Possession of data does not authorize a university to train a model on it or release it to others. Participant permissions, confidentiality, contractual restrictions, and intellectual-property rights may limit particular uses. Agreements among universities, collaborators, and funders should specify who may use the material, for what purposes, who maintains it, what can be published, and how contributors are credited. Some resources can be openly shared; others require controlled access or secure analysis.

\section{Alignment and independent evaluation as research priorities}

University leadership must include investigating how AI behaves and how its development is directed. Research-level capability does not ensure that an AI system pursues the intended purpose or remains open to correction. These questions should not be left primarily to the companies building the most powerful systems. Universities can bring computational, experimental, philosophical, and social-scientific approaches, but meaningful participation requires resources and access to study AI behavior, not opportunities to comment on it from outside.

Reliability, alignment, and governance address different parts of this problem. \emph{Reliability} concerns whether a system performs a task correctly under relevant conditions: whether its calculations and reported measurements are accurate. \emph{Alignment} concerns whether its behavior serves the intended purpose, including when instructions are incomplete or an easily measured target diverges from that purpose. A system rewarded for producing promising findings may favor positive results and neglect contradictory evidence. Such behavior would undermine the scientific purpose of the work even while improving the measured target. \emph{Governance} concerns who sets the purposes and rules, whose interests are considered, and how decisions can be challenged. Accurate calculations resolve neither of the latter two.

Research environments provide concrete settings in which to investigate these issues. Does an AI assistant revise a proposal when evidence contradicts it? Does it disclose uncertainty, or lose important qualifications when passing information between programs? Can researchers inspect, interrupt, and redirect a long sequence of AI-directed tasks? Do incentives to increase reported discoveries encourage shortcuts in checking them? Experiments designed to expose failures, independent replication, and close examination of research records inform these questions.

Independent evaluation requires the ability to report unwelcome findings. Research agreements should protect publication of unfavorable results rather than give a sponsor power to suppress them. Evaluators should disclose relevant financial and institutional relationships, and teams should have authority to investigate failures without depending on those responsible for promoting a successful demonstration. Universities also face funding pressures and incentives to exaggerate progress; their independence must be secured by enforceable arrangements rather than assumed from their academic status. Where consequences extend beyond the laboratory, affected communities should have a voice in oversight.

Access limits what evaluation can establish. Experiments on an openly available model cannot certify the safety or reliability of a different commercial system. Testing through a commercial interface reveals important aspects of that system's behavior, but some questions require access to internal processes or training information. Investigating those questions requires negotiated research access and secure testing arrangements. Sharing research likewise does not mean releasing sensitive data or every potentially harmful capability without safeguards. The aim is a durable capacity for independent investigation and improvement, not the assumption that university ownership makes AI trustworthy.

\section{Organization, allocation, and sustained support}

This program should be governed as research infrastructure, with professional operational leadership. Faculty and researchers from participating disciplines should hold substantial authority over scientific priorities, working with research-computing teams, engineers, librarians, laboratory staff, and specialists in security and data stewardship. Procurement procedures cannot determine the research agenda, and a faculty committee alone cannot operate a dependable service. A consortium should specify who directs its scientific program, who maintains its facilities, who evaluates its performance, and how those responsibilities are reviewed.

Resource allocation should distinguish ordinary participation from unusually intensive work. A basic level of access should be available to all researchers without a project proposal: they need room to explore before they can justify a major investigation. Transparent limits on resource use and appropriate safeguards are compatible with access that is free at the point of use. Projects requiring extensive computation to train AI models, long sequences of AI-directed tool use, large-scale checks of research findings, or substantial laboratory time need a separate allocation process, including a rapid route for promising exploratory work rather than a requirement that every idea arrive as a mature grant proposal.

Allocation should support discovery, checking, and infrastructure development, not only first announcements. Some capacity should be reserved for unfamiliar approaches and for researchers without established access to major resources. Decisions need disciplinary judgment, published criteria, reasons that can be examined, and opportunities for appeal. A consortium restricted in practice to wealthy universities would reproduce the concentration it seeks to reduce. Membership and funding arrangements should create meaningful routes for less-resourced institutions and collaborations with public laboratories and nonprofit research organizations.

A credible budget begins with research needs and existing resources. Running existing models requires computation, storage, networking, security, and continuing engineering support. Developing models adds data preparation, repeated training and evaluation, and research personnel. Connecting AI to laboratories adds instrumentation, integration, maintenance, and safety requirements. Dependable access and the development of new research capabilities can proceed together. The balance should follow scientific priorities, comparative costs, and the requirements for access and preservation identified in Section~4.

An initial budget should specify the intended users, workloads, service period, and research deliverables. Compare purchased services, shared facilities, and local operation using those same requirements. Include staff time, integration, checking and repair, maintenance, and replacement alongside computing charges. Separate one-time development costs from continuing operation. For operational services, report the cost of independently assessed successful tasks as well as the cost of generating output. Research projects should additionally be assessed for the methods, evidence, and capabilities they develop. This gives decisions about expansion a basis in observed use and scientific value.

Personnel are part of the infrastructure. Research engineers, scientific software developers, data stewards, and laboratory specialists need stable career paths and a substantive role in the research they support. Researchers need recognition for evaluations, maintained software, shared datasets, negative findings, and verification work.

As AI changes how research contributions are produced, universities will also need to reconsider how they assess individual contributions and support early-career researchers. Evaluation should recognize the formulation of significant questions, sound scientific judgment, and responsibility for results, alongside the development of methods and shared resources. As AI makes it easier to produce analyses and written reports, rewarding their number becomes an even weaker proxy for scientific value.

Industry partnerships can contribute AI models, computation, expertise, and routes to practical application. Agreements should protect the ability to publish, retain important research records, and continue essential work when a contract ends. Universities should neither reject commercial collaboration nor mistake it for a substitute for their own research capacity. The relevant question is what they remain able to investigate and decide for themselves.

\section{A first program in experimental and mathematical research}

A consortium could begin with dependable baseline access, interdisciplinary development teams, and two complementary research programs. One would connect AI to experimental inquiry; the other would connect mathematical discovery to checking and knowledge organization. These are illustrative starting points; institutions should select programs suited to their research strengths and communities. Each should have a defined scope, scientific and technical leadership, and a review date. Together they would test both directions of the argument: whether AI improves disciplinary research, and whether evidence and problems arising in that research enable university teams to improve AI itself.

In an experimental program, such as materials characterization, a team would connect selected instruments and documented datasets to AI models that help choose informative experiments. Researchers would define permitted procedures, measurement standards, and conditions for human intervention. The program would produce a working integration, documented experimental records, and independently repeated findings. It would assess whether AI improves experimental selection or helps resolve discrepancies between simulation and observation. Failures would inform joint research on predictions under conditions poorly represented in the original data. The CRESt platform described in Section~7 provides a precedent; adaptation to another laboratory would itself be an outcome to assess.

In mathematics, a team could begin with a bounded area in which participating researchers can independently assess the work. It would build a research environment combining a selected collection of mathematical papers and references with AI assistance and proof-checking tools where suitable. Mathematical records would retain links to their original sources. The pilot would produce a usable prototype, an independently assessed set of research tasks, and records of successful and unsuccessful uses. Evaluation would ask whether researchers find relevant prior work, retain the assumptions required when applying or adapting a result, develop useful approaches, and identify consequential errors. Researchers would assess whether formal statements faithfully express the intended mathematics and whether resulting contributions have scientific value. These cases would also support research on mathematical representation, reasoning, and evaluation of AI systems.

Both programs should compare their approach with existing practice, simpler automation, and current commercial AI assistance using matched sources, tasks, and resources where possible. Document the systems and versions used, and refresh comparisons when changes materially affect the question being tested. Evaluation should count total resources consumed and the effort required to check and repair mistakes. Claims about improved AI should identify what changed---the model, its information sources, or its use of tools---and be tested on problems or experimental conditions withheld from development. Independent assessors should examine both the disciplinary findings and the claimed improvements to AI.

Expansion should follow demonstrated value, with criteria for revision stated in advance. A useful specialized model can receive further support; an ineffective service or unreliable laboratory integration should be revised or discontinued. Early programs should leave reusable software and methods, documented findings, and candid assessments of their limitations, shared where rights and confidentiality permit. Evidence that another group can use these resources successfully would strengthen the case for further investment.

\section{From access to leadership}

Changes already visible in mathematics give urgency to preparation across academia. Universities can build scientific and technical capacity now, revising their commitments as research establishes what works.

University leaders can begin by assessing existing computing and laboratory capacity, joining or forming research consortia, funding interdisciplinary AI-development teams, and establishing dependable baseline access. Scientific leadership requires authority over these commitments and budgets that sustain their operation. An institutional strategy should identify the research capabilities it intends to build, the people responsible for developing them, and the evidence on which continued investment will depend.

Funding agencies should support shared computation, foundational AI research, experimental facilities, knowledge infrastructure, and independent evaluation as connected parts of the research system. National centers can supply resources beyond individual institutions while universities contribute scientific direction and development expertise. Funding should sustain continuing operations and maintenance as well as initial demonstrations, allowing promising work to expand and unsuccessful designs to change.

Researchers, disciplinary societies, and academic consortia can articulate requirements, coordinate shared resources, and represent the research community in decisions about access and standards. Collaboration should extend beyond a few wealthy institutions and beyond academia itself. Leadership means making consequential contributions and accepting responsibility for them, not claiming exclusive authority over science or AI.

The ambition is to bring more of the world's intellectual capacity into the development of AI and the advancement of knowledge: more questions, more methods, more independent scrutiny, and more communities able to contribute. Universities should sustain the infrastructure and professional communities through which discoveries remain available for others to examine and develop.

\section*{AI assistance disclosure}

This article was developed with assistance from ChatGPT (OpenAI) and Claude (Anthropic). These tools were used to revise text, to test arguments and supply counterarguments, and to help locate and check sources. The author determined the positions advanced here, made the final decisions on all text, and is responsible for the content.

\end{document}